\documentclass[twocolumn,10pt]{IEEEtran}
\usepackage{ifpdf, flushend,subfigure,tabularx}

\ifCLASSINFOpdf
  \usepackage[pdftex]{graphicx}
  \graphicspath{{../pdf/}{../jpeg/}}
  \DeclareGraphicsExtensions{.pdf,.jpeg,.png}
\else
  \usepackage[dvips]{graphicx}
  \graphicspath{{../eps/}}
  \DeclareGraphicsExtensions{.eps}
\fi
\usepackage{graphicx}
\usepackage{epstopdf}
\graphicspath{{../}}
\usepackage{multirow}
\usepackage{float}
\usepackage[cmex10]{amsmath}
\usepackage {amssymb}
\usepackage{array}
\usepackage{mdwmath}
\usepackage{mdwtab}
\usepackage{eqparbox}
\usepackage{nomencl}
\usepackage{cite}
\usepackage{algorithm}
\usepackage{algorithmic}
\usepackage{makeidx}
\usepackage{ifthen}
\usepackage{subfigure}
\usepackage{caption}
\usepackage{pdflscape}
\usepackage{hyperref}
\usepackage{xcolor}

\def\bg{{\mathbf{g}}}
\def\bx{{\mathbf{x}}}
\def\bz{{\mathbf{z}}}

\def\bu{{\mathbf{u}}}
\def\bq{{\mathbf{q}}}

\newcommand{\argmin}{\operatornamewithlimits{argmin}}

\newcommand{\beq}{\begin{equation}}
\newcommand{\eeq}{\end{equation}}
\newcommand{\beqn}{\begin{eqnarray}}
\newcommand{\eeqn}{\end{eqnarray}}
\newcommand{\beqno}{\begin{eqnarray*}}
\newcommand{\eeqno}{\end{eqnarray*}}
\newcommand{\bma}{\begin{displaymath}}
\newcommand{\ema}{\end{displaymath}}
\newcommand{\bnu}{\begin{enumerate}}
\newcommand{\enu}{\end{enumerate}}
\newcommand{\bce}{\begin{center}}
\newcommand{\ece}{\end{center}}
\newcommand{\btb}{\begin{tabular}}
\newcommand{\etb}{\end{tabular}}

\begin{document}


\title{A Fairness-Aware Attacker-Defender Model for Optimal Edge Network Operation and Protection}

\author{\IEEEauthorblockN{Duong~Thuy~Anh~Nguyen,~\IEEEmembership{Student Member,~IEEE}, Jiaming~Cheng,~\IEEEmembership{Student Member,~IEEE}, ~Ni~Trieu,~\IEEEmembership{Member,~IEEE}, and~Duong~Tung~Nguyen,~\IEEEmembership{Member,~IEEE}}  
\thanks{Duong Thuy Anh Nguyen, Ni Trieu, and Duong Tung Nguyen are with 
Arizona State University, Tempe, AZ, United States. Email: \textit{\{dtnguy52,~nitrieu,~duongnt\}@asu.edu.}
Jiaming Cheng is with 
University of British Columbia, Vancouver, BC, Canada. Email: \textit{jiaming@ece.ubc.ca}.
 }
 }

 

\maketitle

\begin{abstract}

While various aspects of edge computing  (EC) have been studied extensively, the current literature has overlooked the robust edge network operations and planning  problem. 
To this end, this letter proposes a novel fairness-aware attacker-defender model for optimal  edge network operation and hardening against possible attacks and disruptions.
The proposed model helps 
 EC platforms identify the set of most critical  nodes to be protected  to mitigate the impact of failures on system performance.
Numerical results demonstrate that the proposed solution not only ensures  good service quality but also maintains   fairness among different areas during 
disruptions. 

\end{abstract}
\begin{IEEEkeywords}
Edge computing,  attacker-defender, fairness, network hardening, proactive protection, network disruptions.
\end{IEEEkeywords}

\printnomenclature


\allowdisplaybreaks
\section{Introduction}
\label{Sec:Intro}

Edge computing (EC) has been proposed as a complement to the cloud to enhance user experience and support a variety of Internet of Things (IoT) applications. 
The new EC paradigm offers powerful computing resources by enabling computation to be performed at the network edge in close proximity  to  end-users, instead of being centralized in cloud environments. The decentralization of resources has alleviated  bandwidth constraints and overcomes latency issues, making the optimization of edge network operation and planning a critial problem. 
%
%
In EC, there are multiple heterogeneous edge nodes (ENs) with varying capacities, thus, an important question is to allocate edge resources efficiently and effectively \cite{Nguyen2018a, Cheng2021, tara,Nguyen2022}, with the aim of reducing latency, energy consumption, and unmet demand. 

In light of the open nature of the EC ecosystem and the vulnerability of IoT devices with limited resources to attacks, reliability and resilience are also crucial challenges, 
 given the growing risk of  failures. 
In practice, various  factors such as cyber-attacks,  natural disasters, software errors, and power outages can result in component failures and  disrupt the edge network operation, leading to substantial economic losses and safety concerns for mission-critical services, and negatively impacting the user experience. 
Despite this, the issue of resilience in EC remains unresolved \cite{Roman2018}.

An initial discussion of reliability challenges in EC was introduced in \cite{Madsen2013}; however, limited efforts have been made to address these challenges. 
One line of research aims at accurately estimating  failures in the field. For instance, in \cite{Aral2017}, a Bayesian-based technique is introduced to exploit the causal relationship between different types of failures and identify the availability of virtual machines running on different data centers. 
In \cite{Aral2021}, the authors leverage spatio-temporal failure dependencies among edge servers and employ a  learning approach to compute the joint failure probability of the servers. Additionally, in \cite{FL_com1, FL_com2}, the authors present a federated learning-based decentralized communication approach  to mitigate the negative impacts on the network after node failures.

Another line of research focuses on optimizing edge systems under unpredictable failures. Reference \cite{Ford2017} proposes to model latency and resilience requirements of EC applications as an optimization problem by reserving bandwidth and resource capacities. To investigate optimal service placement strategies,  \cite{Qu2021,Cheng2021} formulate robust optimization (RO) models under uncertain EN failures. A heuristic solution is proposed in \cite{Ford2017}. The later work \cite{Qu2021} exploits the monotone submodular property to develop approximation algorithms while a globally optimal solution is proposed in \cite{Cheng2021}. 
\textbf{Reference \cite{fhe20} presents a \textit{min-max} optimization model for preventive priority setting to handle load balancing against  controller failures in software-defined networks.}
A robust mixed integer linear programming (MILP) model is proposed in \cite{mito211}  to minimize the total required backup capacity against simultaneous failures of physical machines.

The existing works primarily focus on the probability of failures and optimizing decisions under failure events using probabilistic or robust models. There is  limited attention paid to developing effective defense strategies prior to failures that can enable the edge platform to dynamically adapt its operations and maintain a high quality of service during disruptions. 
this letter seeks to address the gap in the literature by focusing on the proactive edge network protection problem to mitigate the impacts of failure events on system performance. 

The protection of edge networks against attacks is a complex issue, due to the uncertainty of network failures. 
For simplicity, 
this work considers EN failures only. 
Node failures can occur due to various reasons including power outages, internal component failures, software errors or misconfigurations, natural disasters, or cyberattacks. When an EN fails, 
it is considered as being attacked by a hypothetical attacker and is virtually eliminated from the network. 
As a result, the workload previously assigned to the failed EN must be reallocated to other active ENs, potentially leading to increased latency for some users. Additionally, as the number of node failures increases, the available edge resources decrease, potentially resulting in an increase in unmet demand. Therefore, it is of utmost importance to mitigate the degradation of services in the face of network failures.

 


\textbf{Contributions:} 
This letter presents a novel bilevel attacker-defender optimization model and corresponding solution techniques for ensuring the survivability of edge networks (ENs) in the face of failures through proactive protection. The objective of the attacker is to maximize the degradation of user experience, as reflected by unmet demand and delay. By considering the problem from the attacker's perspective, the proposed techniques help to understand the most disruptive attacks and identify the critical set of ENs for effective protection.


Protection of ENs can be achieved through various means, such as placement in secure locations, installation of firewalls and security software, implementation of intrusion detection and prevention systems, and use of Uninterruptible Power Supply (UPS) units. The proposed approach solves a max-min model representing the optimal attack problem to obtain a robust proactive protection strategy, using an efficient algorithm based on linear programming duality to compute the exact global optimum of the bi-level optimization problem.


In addition, node failures can lead to insufficient edge resources and raise concerns about the fair distribution of resources to users during disruptions. The issue of fairness-aware resource allocation under disruption has been overlooked in prior work. To address this, a novel notion of fairness in the context of edge resource allocation is introduced,
which prevents the edge platform from prioritizing serving the demands from some areas over others to maximize its utility and ensures
sufficient resources to serve critical services in every area
during network failures. 
The fairness aspect of the decision-making process is reflected through the fairness constraints in our formulation, which allows the platform to balance the trade-off between fairness and efficiency. Simulation results show that proactive provisioning significantly improves  system performance in terms of service quality and fairness. 

\section{System Model and Problem Formulation}
\label{Sec:ProblemFormulation}

\subsection{System Model}
\begin{figure}[hbt!]
	\centering
       \includegraphics[width=0.40\textwidth,height=0.11\textheight]{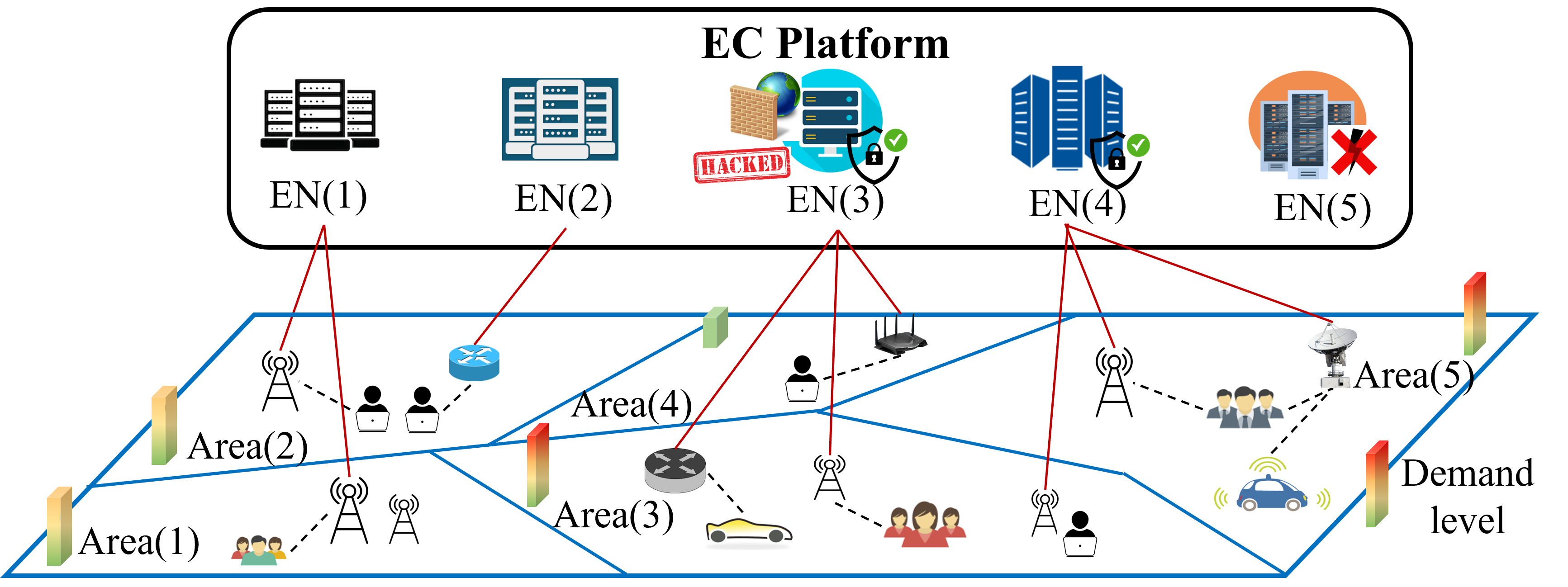}
         	\caption{System model}
	\label{fig:model}
\end{figure}

The system model is illustrated in Fig.~\ref{fig:model}.
We consider an  EC platform that manages a set $\mathcal{N}$ of $N$ ENs to serve users located in different areas, each of which 
is represented by an Access Point (AP). Let $\mathcal{M}$ be the set of $M$ APs.  
We define $i$ and $j$ as the AP index and EN index, respectively. The computing resource capacity of EN $j$ is $C_j$. Let  $\lambda_i$ be the total resource demand of users in area $i$. We use the terms ``\textit{demand}" and ``\textit{workload}" interchangeably in this letter. The network delay between AP $i$ and EN $j$ is $d_{i,j}$. The platform aims to reduce the overall delay as well as maximize the amount of workload served at the ENs. 
Denote the workload from area $i$ that is allocated to EN $j$ by $x_{i,j}$. The service requests from each area must be either served by some ENs or dropped (i.e., counted as the unmet demand $q_{i}$), and the penalty is $\phi_i$ for each unit of unmet demand in area $i$. 

\vspace{-0.1cm} 

The operation of the edge network can be disrupted for various reasons. 
To account for this, we consider the presence of a hypothetical intelligent adversary, i.e., an attacker, whose objective is to cause maximum damage to the edge system performance given its limited attack resources.
The EC platform acts as the defender. We assume the attacker has complete knowledge of the  optimal operation plan of the platform. This is a reasonable assumption as the defender would not be worse off even if the attacker possesses a less-than-perfect model of the defender's system or fails to execute a perfect attack plan. 

When an EN fails due to any reason, 
we consider it as being attacked by the hypothetical attacker. 
The workload initially assigned to the failed EN must  be reallocated to the remaining active ENs. The survivability of the system against edge failures can be achieved through EN protection.  We assume that the EC platform has a certain budget for EN hardening, and the maximum number of ENs that can be protected is represented by $K$. 
Then, the platform can optimize the EN protection decision by viewing this problem through the lens of an attacker who 
can only successfully attack at most $K$ simultaneous ENs. 
Obviously,  $K\leq |\mathcal{N}|=N$.

Since the exact node failures can not be predicted accurately at the time of making the protection decision, we propose techniques to search for 
 optimal disruptive attacks, given the posited offensive resources of the attacker. To this end, we formulate a bilevel optimization  problem where the defender solves the inner  problem to minimize the impact on user experience, while the attacker solves the outer problem to maximize the  degradation of system performance. The solution to this problem yields the most disruptive attack that the attacker could undertake, as well as  the  critical ENs that should be reinforced for optimal protection. 

\vspace{-0.1cm}
\subsection{EC Platform Optimization Model}
First, we present the optimal operating model for the EC  platform under normal conditions without attack. The platform's objective is to improve the  user experience by reducing the overall delay and the extent of unmet demand. 
Hence, it needs to solve the following optimization problem:
\begin{subequations}
\begin{align}
\label{eq-Obj:D-Prob}
\textbf{(D)}~~ \min_{\bx,\bq}~ & (1-\gamma)\sum_{i}\phi_i q_{i} + \gamma\sum_{i,j} d_{i,j} x_{i,j} \\
\text{s.t.}~\label{eq-const:Capacity}
& \sum_{i}x_{i,j}\leq C_j,~\forall j\\
\label{eq-const:WorkloadAllocation}
& \sum_{j}x_{i,j}+q_{i}= \lambda_i,~ \forall i; ~x_{i,j} \leq C_ja_{i,j}   ,~ \forall i,j\\
\label{eq-const:fair}
& \frac{q_{i}}{\lambda_i} \leq \theta,~ \forall i;~~ \Bigg| \frac{q_{i}}{\lambda_i} - \frac{q_{i'}}{\lambda_{i'}} \Bigg| \leq \beta,~ \forall i,i'\\
& x_{i,j} \geq 0,~\forall i,j, ~~ q_{i} \geq 0,~\forall i. \label{eq-const:vars}
\end{align}
\end{subequations}
The first term in the objective function (\ref{eq-Obj:D-Prob}) captures the total penalty for unmet demand. 
The second term represents the  delay penalty. The weighing factor $\gamma$ can be used by the EC platform to express its priorities towards reducing the delay and the unmet demand. Overall,  the EC platform aims to minimize the weighted sum of the unmet demand penalty and the delay penalty. The computing resource capacity constraints are presented in \eqref{eq-const:Capacity}, which shows that the total amount of workload allocated to EN $j$ can not exceed its capacity. 
The workload allocation constraints \eqref{eq-const:WorkloadAllocation} indicate the demand from each area must be either served by some ENs or dropped. 

Since we focus on delay-sensitive services, constraints  \eqref{eq-const:WorkloadAllocation} also impose that the demand from each area can only be served by a subset of ENs which are sufficiently close to that area. We use $a_{i,j}$ as a binary indicator to denote if EN $j$ satisfies the prerequisites for serving the demand from area $i$. In particular, if $a_{i,j} = 0$, then $x_{i,j} = 0$, i.e., EN $j$ is not eligible to serve the demand from area $i$.  Otherwise, the workload assigned from area $i$ to EN $j$ should not exceed the capacity $C_j$ of the EN.

For quality control, the EC platform can enforce that the ratio of the unmet demand to the total demand in each area $i$ should not exceed a given threshold $\theta$ as shown in 
\eqref{eq-const:fair}. This condition ensures that each area is allocated sufficient resources for operating critical services in case of disruption. Constraints \eqref{eq-const:fair} enforces a soft fairness condition among different areas in terms of the proportions of unmet demand. Specifically, the proportions of unmet demand between any two areas should not differ more than a certain threshold $\beta$. 

\vspace{-0.2cm}
\subsection{Attacker-Defender Model}
We are now ready to describe the attacker-defender model to determine the set of critical ENs. The attacker can use this model to find an optimal attack plan, given its limited offensive resources, to maximally degrades the performance of the EC platform during disruptions by increasing the total amount of unmet demand as well as the overall network delay
By studying attack strategies through the lens of the attacker, it helps us understand how to make the edge network less vulnerable. 
In our model,  the EC platform acts as the defender.

Let $z_j$ be a binary variable that takes the value of $1$ if EN $j$ is successfully attacked by the attacker. For simplicity, we assume that an EN will completely fail if it is attacked. 
Denote an attack plan 
of the attacker by $\bz=(z_1,\cdots,z_N)$ $\in \mathcal{Z}$, where $\mathcal{Z}$ is the set of all possible attack plans that can be carried out by the attacker. In other words, $\mathcal{Z}$ is the set of the attacker's options to interdict different combinations of ENs in the network. A natural form of set $\mathcal{Z}$ is: 
 \begin{align}\label{eq-set:AttackerSet}
    \mathcal{Z}=\Big\{\bz\in\{0,1\}^{N}:~\sum_j z_j\leq K\Big\},
\end{align}
where $K$ is the maximum number of ENs that the attacker can attack simultaneously. 
Note that  our proposed model and solution can easily deal with a more general (polyhedral) form of set  $\mathcal{Z}$. 
Overall, the attacker-defender model can be expressed as a max-min optimization problem as follows:
\beqn
\label{ADProb}
\textbf{(AD)}~~ \max_{\bz \in \mathcal{Z}} ~\min_{\bx,\bq\in\mathcal{S}(\bz)} ~ (1-\gamma)\sum_{i}\phi_i q_{i} +\gamma\sum_{i,j} d_{i,j} x_{i,j},
\eeqn
where  $\mathcal{S}(\bz)$ represents the defender's set of feasible actions, restricted by the attack plan $\bz$. Concretely, 
 $\mathcal{S}(\bz)$ captures all the operational constraints of the defender during an attack. Since attacked ENs cannot serve any demand, constraints (\ref{eq-const:Capacity}) 
need to be modified as follows:
\beqn
\label{eq-const:CapacityA}
\sum_{i}x_{i,j}\leq C_j(1-z_j),\quad \forall j.
\eeqn
It can be seen that if EN $j$ fails  (i.e., $z_j=1$),  then \eqref{eq-const:CapacityA} implies  $\sum_{i}x_{i,j} = 0$, thus, $x_{i,j} = 0,~\forall i,j$, and EN $j$ is unable to serve any user requests. Since  \eqref{eq-const:CapacityA} already enforces  $x_{i,j} = 0,~\forall i,j$, when $z_j=1$, we do not need to modify  (\ref{eq-const:WorkloadAllocation}). Hence, the feasible set $\mathcal{S}(\bz)$ given  the attack plan $\bz$ can be given as:
\beqn
\label{eq-set:S(z)}
\mathcal{S}(\bz)=\Big\{\bx,\bq:~\eqref{eq-const:WorkloadAllocation} - \eqref{eq-const:vars}, ~ \eqref{eq-const:CapacityA} \Big\}.
\eeqn
The attacker-defender model \eqref{ADProb} aims to find a critical set of ENs by identifying a maximally disruptive resource-constrained attack that an attacker might undertake. The inner minimization problem in \eqref{ADProb} represents the optimal operation model of the EC platform (i.e., the defender) under a specific attack plan $\bz$. The outer maximization problem models the attacker's objective to incur the highest cost to the EC platform during the attack. 
The solutions obtained from solving the optimal-attack problem \eqref{ADProb} suggest the set of $K$-critical ENs to be protected to maximize the user experience.

The attacker-defender problem \textbf{(AD)} in  \eqref{ADProb}  is a challenging bilevel program \cite{bilevel1}. Refer to \textit{Appendix A} in \cite{report} for the bilevel formulation. The computational difficulty stems from the \textit{max-min} structure of the problem. The outer problem also referred to as the leader problem, represents the attacker problem that seeks to maximize the defender's loss caused by the disruption. The inner minimization problem, also referred to as the follower problem, captures the optimal operation of the edge platform under a specific attack plan $\bz$. The platform aims to optimally allocate the workload to non-attacked ENs to mitigate the impact of attacks by minimizing the amount of unmet demand as well as the overall network delay. 

\allowdisplaybreaks
\section{Solution Approach}
\label{Sec:Algorithm}
To solve 
the problem \eqref{ADProb}, 
 we first convert the inner  problem into a linear program (LP) for a fixed value of $\bz$. 
Then, \eqref{ADProb} becomes a bilevel LP with binary variables in the upper level, which imposes difficulties \cite{bilevel1}. 
We propose transforming the inner minimization problem into a maximization problem by relying on  LP duality. Consequently, the \textit{max-min} problem \eqref{ADProb} can now be reformulated as a \textit{max-max} problem, which is simply a maximization problem. The resulting maximization problem is a mixed-integer nonlinear program (MINLP) containing several bilinear terms. By linearizing the bilinear terms, we can convert the MINLP into an MILP that can be solved  by off-the-shelf solvers. The detailed solution is described below.


First, given a fixed attack plan $\bz$, 
the inner minimization problem in \eqref{ADProb} can be rewritten as an equivalent LP below: 
\begin{subequations}
\begin{align}
\label{eq:InnerMinimization}
\min_{\bx,\bq} \quad   &(1-\gamma)\sum_{i}\phi_i q_{i} + \gamma\sum_{i,j}d_{i,j}x_{i,j} &&  ~~\\
\text{s.t.} \quad 
&\sum_{i}x_{i,j}\leq C_j(1-z_j), ~\forall j  \label{eq-const:SP1} ~  &&(\pi_j)\\
\label{eq-const:SP2}
&\sum_{j}x_{i,j}+q_{i} = \lambda_i,~ \forall i ~  &&(\mu_i)\\
&x_{i,j} \leq C_ja_{i,j}, ~\forall i,j \label{eq-const:SP3}   &&(\sigma_{i,j})\\
&\frac{q_{i}}{\lambda_i} - \frac{q_{i'}}{\lambda_{i'}} \leq \beta,~ \forall i,i'    &&(\eta_{i,i'})\label{eq-const:SP4}\\
&\frac{q_{i}}{\lambda_i} - \frac{q_{i'}}{\lambda_{i'}}  \geq -\beta,~ \forall i,i'   &&(\tau_{i,i'}) \label{eq-const:SP5}\\
&x_{i,j} \geq 0,\forall i,j;~ q_{i} \geq 0,\forall i;~\frac{q_{i}}{\lambda_i}  \leq \theta, \forall i       &&(\nu_{i})\label{eq-const:SP6}
\end{align}
\end{subequations}
$\pi_j,\mu_i,\sigma_{i,j},\eta_{i,i'},\tau_{i,i'},\nu_{i}$ are the dual variables associated with constraints (\ref{eq-const:SP1})-(\ref{eq-const:SP6}), respectively. 
To transform the \textit{max-min} bilevel program into a single-level  problem, a common approach is to utilize the  Karush–Kuhn–Tucker (KKT) conditions to convert the inner follower problem into a set of linear constraints (see \textit{Appendix B} in our technical report \cite{report}). 
However, this solution is not computationally efficient because the reformulated single-level optimization problem contains a large number of complimentary constraints. Instead, we    employ  the  LP duality to convert the minimization problem (\ref{eq:InnerMinimization})-(\ref{eq-const:SP6}) to the following equivalent maximization problem.
\begin{subequations}
\begin{align}
\label{eq:DualityReformulation}
\max_{\pi,\sigma,\mu,\eta,\tau,\nu}  -\!\sum_{j}\!C_j(1\!-\!z_j)\pi_j\!+\!\sum_{i}\!\lambda_i\mu_i\!-\!\sum_{i,j}\!C_ja_{i,j}\sigma_{i,j} \nonumber \\ 
-\sum_{i}\sum_{i'=i+1}^M\beta(\eta_{i,i'}+\tau_{i,i'})-\sum_{i}\theta\nu_{i} \\
\text{s.t.} \quad  \quad \quad  \quad  \quad \quad  - \pi_j + \mu_i - \sigma_{i,j} \le \gamma d_{i,j}, ~\forall i,j    \label{eq-const:Duality1}\\
\label{eq-const:Duality2}
\mu_i-\frac{1}{\lambda_i}\sum_{l=i+1}^M\left(\eta_{i,l}-\tau_{i,l}\right) +\frac{1}{\lambda_i}\sum_{l=1}^{i-1}\left(\eta_{l,i}-\tau_{l,i}\right) -\dfrac{\nu_{i}}{\lambda_i}  \nonumber \\ 
\le (1-\gamma)\phi_i,~ \forall i\\
\pi_j \geq 0,\sigma_{i,j} \geq 0,\eta_{i,i'} \geq 0,\tau_{i,i'} \geq 0,\nu_{i} \geq 0,~\forall i,i',j.\label{eq-const:Dualitye}
\end{align}
\end{subequations}
The bilinear term $(1-z_j)\pi_j$ in \eqref{eq:DualityReformulation} is a product of a continuous variable and a binary variable. To tackle this, we define a new auxiliary  non-negative continuous variable $g_j$, where $g_j=(1-z_j)\pi_j,~ \forall j$. The constraint $g_j=(1-z_j)\pi_j$ can be  implemented equivalently through the following linear inequalities:
\begin{subequations}
\begin{align}
g_j \le M_j (1-z_j), ~~ g_j\le\pi_j &, ~\forall j\label{eq-const:Duality3}\\
\label{eq-const:Duality4}
g_j \ge 0, ~~ g_j\ge \pi_j - M_j z_j &,~ \forall j,
\end{align}
\end{subequations}
where each $M_j$ is a sufficiently large number. Based on the linearization steps above, we obtain an MILP that is equivalent to the attacker-defender problem (\ref{ADProb}), as follows:
\begin{subequations}
\begin{align}
\label{eq:MILPDualObj}
\max_{\bz,\bg,\pi,\sigma,\mu,\eta,\tau,\nu}  &\!\!\!-\!\sum_{j}\!C_jg_j\!+\!\sum_{i}\!\lambda_i\mu_i\!-\!\sum_{i,j}\!C_ja_{i,j}\sigma_{i,j}\cr
&\!\!\!\!\!\!-\sum_{i}\sum_{i'=i+1}^M\beta(\eta_{i,i'}+\tau_{i,i'})-\sum_{i}\theta\nu_{i}\\
\text{s.t.} \quad\quad 
        &\!\!\!\eqref{eq-const:Duality1}-\eqref{eq-const:Duality4};\sum_j z_j\leq K;z_j\in\{0,1\},\forall j. \label{eq:MILPDualConst}
\end{align}
\end{subequations}

This MILP problem can now be solved efficiently by MILP solvers  such as Gurobi, CPLEX, and Mosek.
Table~\ref{table:SizeComparision} compares the number of constraints and variables between duality- and KKT- based reformulations. The duality reformulation has a simpler structure and fewer variables and constraints and is usually less computationally expensive. Another disadvantage of the KKT-based method lies in  the weak relaxation of the big-M method, which can greatly affect the running time. 

\begin{table}[h]
\centering
\begin{tabular}{|c|c|c|}
\hline
                    & KKT           & Duality     \\ \hline
\# constraints      & $5N+M(8M+8N+1)$ & $6N+2M(M+N)$ \\ \hline
\# variables ($\mathbb{Z}$)       & $2N+2M(M+N)$  & $N$          \\ \hline
\# variables ($\mathbb{R}$) & $N+2M(M+N+1)$ & $2N+M(2M+N)$ \\ \hline
\end{tabular}
\caption{KKT-based vs duality-based reformulations}
\label{table:SizeComparision}
\vspace{-0.2cm}
\end{table}

\section{Numerical Results}
\label{Sec:NumericalResults}
Similar to  
\cite{Nguyen2018a,Nguyen2022,tara}, we adopt the  Barabasi-Albert model 
to generate a random scale-free edge network topology with $100$ nodes and the attachment rate of $2$ \cite{Nguyen2018a}. The link delay between each pair of nodes is  generated randomly in the range of [2, 5] ms. The network delay between any two nodes is  the delay of the shortest path between them \cite{tara}. From the generated topology, 80 nodes are chosen as APs and 30 nodes are chosen as ENs for our performance evaluation (i.e., $M=80$, $N=30$). 
Also, $a_{i,j}$ is set to be 1 if $d_{i,j}$ is less than $20$ ms. 
The  capacity $C_j$ of each EN $j$ is randomly selected from 
the set $\{16,32,64,128,256,512,1024\}$ vCPUs. 
The resource demand (i.e., workload)
in each area is randomly drawn in the range of [20, 35] vCPUs \cite{tara}.
The penalty for each unit of unmet demand  is set to be 5 \$/vCPU (i.e.,  $\phi_i = 5,~\forall i$).
The EC platform can adjust 
$\gamma$ 
to control the trade-off between the unmet demand  and the  overall network delay.
A smaller value of $\gamma$ implies the platform gives more priority to unmet demand mitigation than to delay reduction.
For illustration purposes, we set $\gamma$ to be $0.1$. 
Also, $\theta$ is set to be 0.8. 

We first evaluate the performance of  the proposed model and
three benchmark schemes, including: (i) \textit{No hardening}: proactive EN protection is not considered; (ii) \textit{Heuristic}: the system protects $K$ ENs with highest computing capacities; and (iii) \textit{Random:} a random protection plan within the protection budget, i.e., protecting $K$ ENs at random. 
The simulation is conducted for different values of $K$ and the  number of \textit{actual} EN failures. 
We run each scheme to obtain a set of ENs to be protected. For our proposed model, we solve 
(\ref{ADProb}) 
to obtain the set of $K$ critical ENs.  This set and the set of actual EN failures are then 
used as an input to problem (\ref{eq:InnerMinimization})-(\ref{eq-const:SP6}) to compute the actual total cost of the platform. For each experiment, we  randomly generate $500$ scenarios of actual EN failures and compute both average and worst-case performance. 
\begin{figure}[h!]
\vspace{-0.4cm}
		\subfigure[Average Cost]{ \includegraphics[width=0.24\textwidth,height=0.10\textheight]{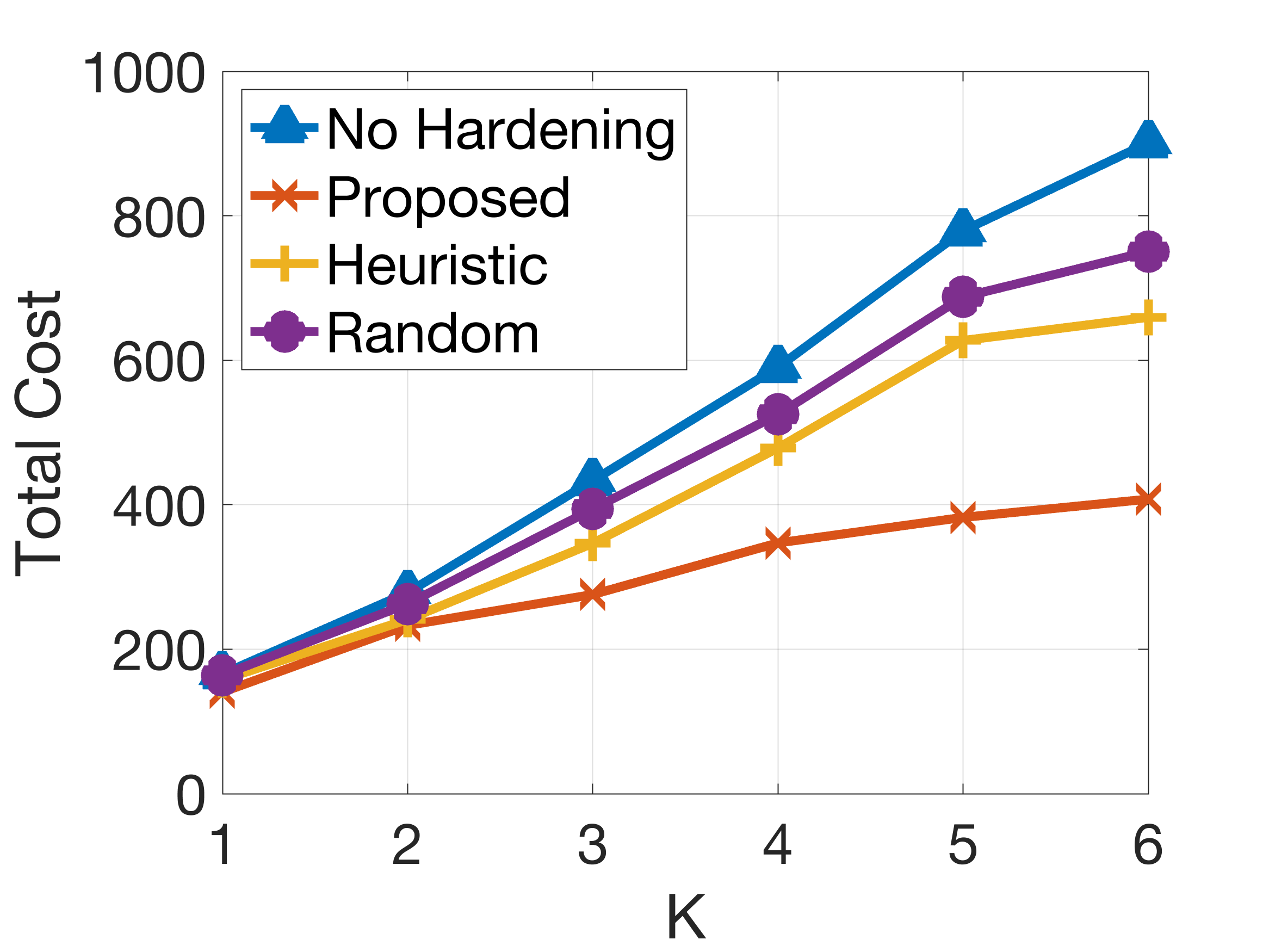}
	    \label{fig:Compareavg}
	}   \hspace*{-2.1em} 
		 \subfigure[Worst-case Cost]{
	     \includegraphics[width=0.24\textwidth,height=0.10\textheight]{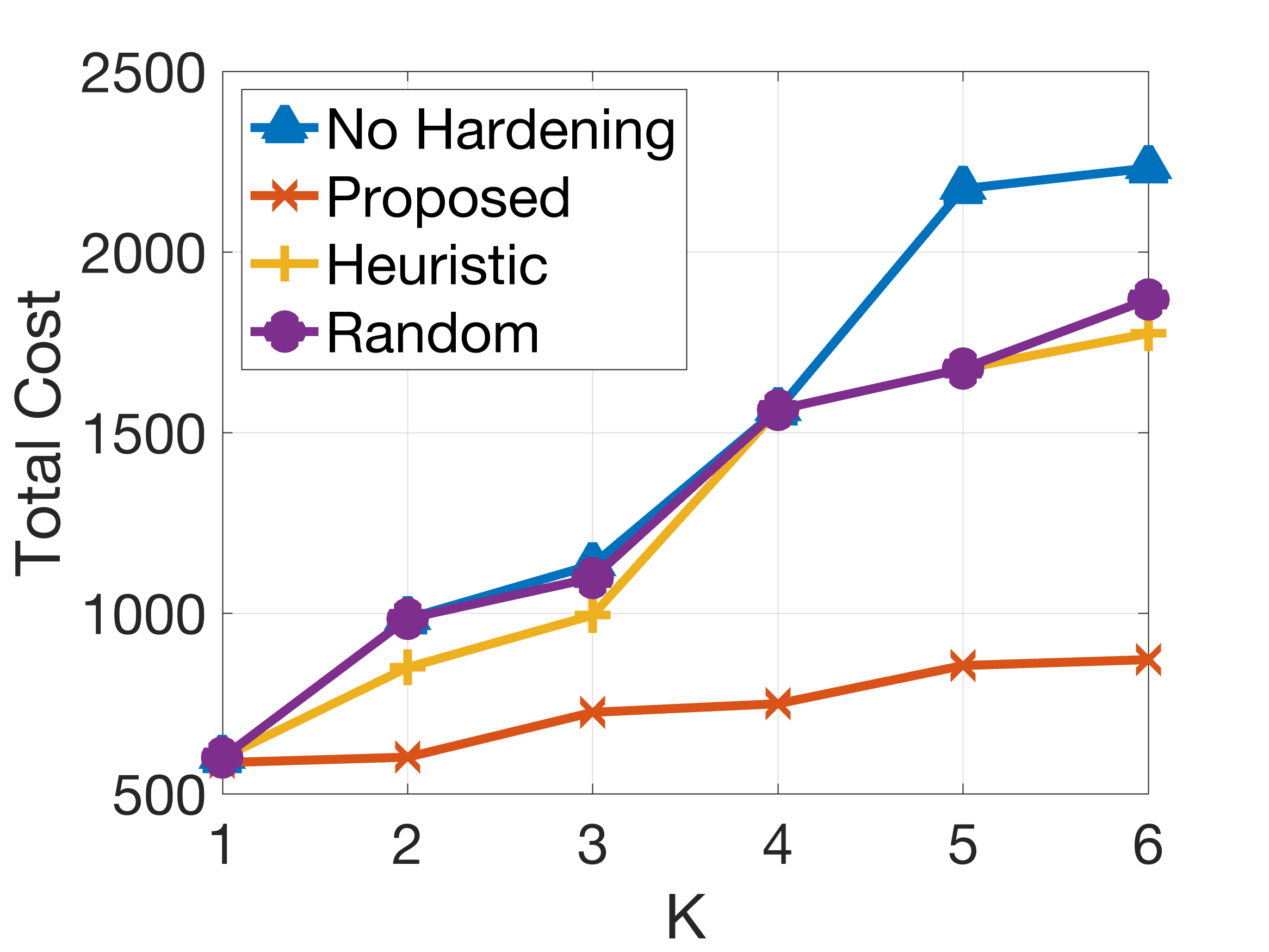}
	     \label{fig:Compareworst}
	}  \vspace{-0.2cm}
	\caption{Performance comparison} 
	\label{fig:4plans}
	\vspace{-0.2cm}
\end{figure}

Fig.~\ref{fig:4plans} compares the performance of the four hardening schemes when the actual number of  failures equals  the defense budget $K$. 
Among the total of $N=30$ ENs, we vary the number of EN failures $K$ from 1 to 6. 
As expected,    ``\textit{No Hardening}''  results in the highest cost for every value of $K$, as all ENs are susceptible to attacks and the attacker has more flexibility to target the critical ENs that cause the greatest disruption. In contrast, 
the attacker has more limited  options with EN protection considered in the other schemes. 
The proposed scheme significantly outperforms the other schemes. 
Also, for every scheme,  the total cost increases as $K$ increases since there are fewer surviving ENs  to serve the demand. 



\begin{figure}[h!]
\vspace{-0.4cm}
		\subfigure[Average Cost]{
		  \includegraphics[width=0.24\textwidth,height=0.09\textheight]{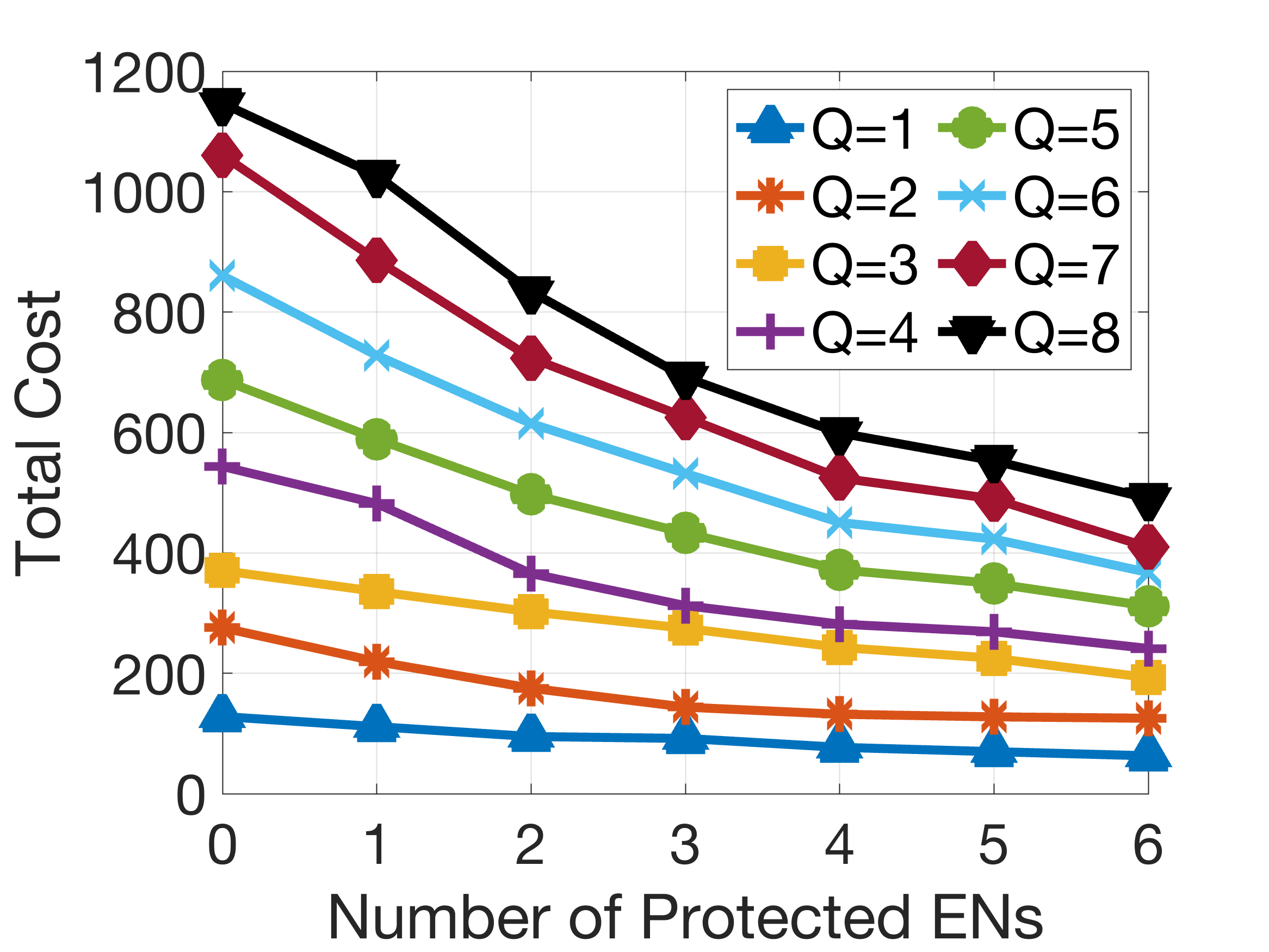}
	    \label{fig:CompareAvgProposed}
	}   \hspace*{-2.1em} 
		 \subfigure[Worst-case Cost]{
	     \includegraphics[width=0.24\textwidth,height=0.09\textheight]{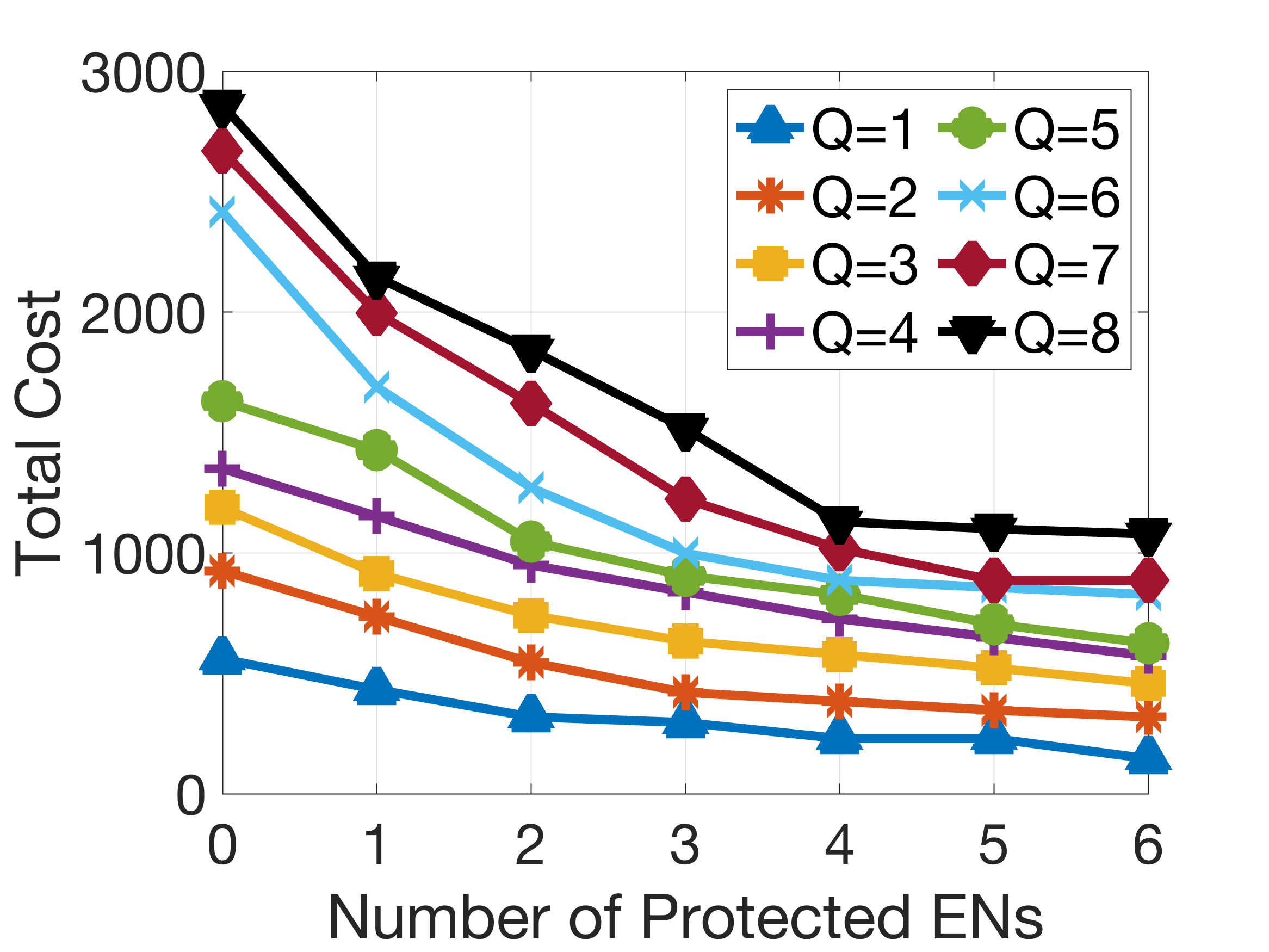}
	     \label{fig:CompareWorstProposed}
	}  \vspace{-0.2cm}
	\caption{Impacts of proactive EN protection on the total cost}
	\label{fig:CompareProtected}
	\vspace{-0.2cm}
\end{figure}

Fig.~\ref{fig:CompareProtected} illustrates the benefit of proactive EN hardening using our model
by showing the effect of having different numbers of protected ENs during the planning stage for different  attack/failure scenarios. 
The performance for different numbers of protected ENs ($K$) and the actual number of EN failures ($Q$) are examined. 
Each curve in the figure corresponds to a given value of $Q$. 
We can easily see that the total cost increases as $Q$ increases. More interestingly,  the defender can drastically reduce its cost by protecting just  a few ENs. For example, protecting only one or two ENs results in a significantly lower total cost compared to the case without any protection. Thus, \textit{proactive EN protection is beneficial even when the number of actual failures is larger than the defense budget}. 

\begin{figure}[h!]
\vspace{-0.4cm}
		\subfigure[$\beta = 0.2$]{
		  \includegraphics[width=0.24\textwidth,height=0.08\textheight]{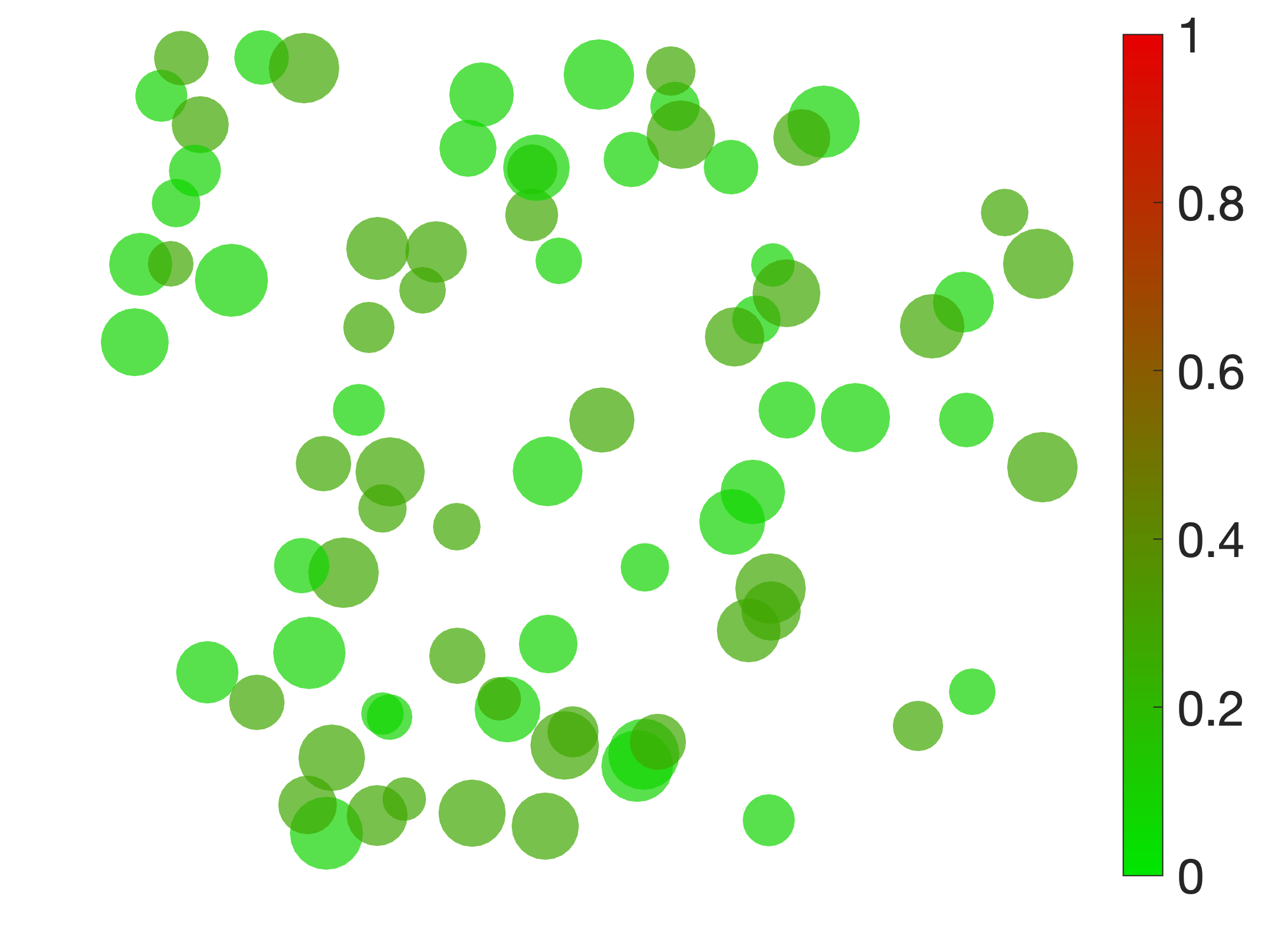}
	    \label{fig:CompareUnmet0.2}
	}   \hspace*{-1.8em} 
		 \subfigure[$\beta = 0.8$]{
	     \includegraphics[width=0.24\textwidth,height=0.08\textheight]{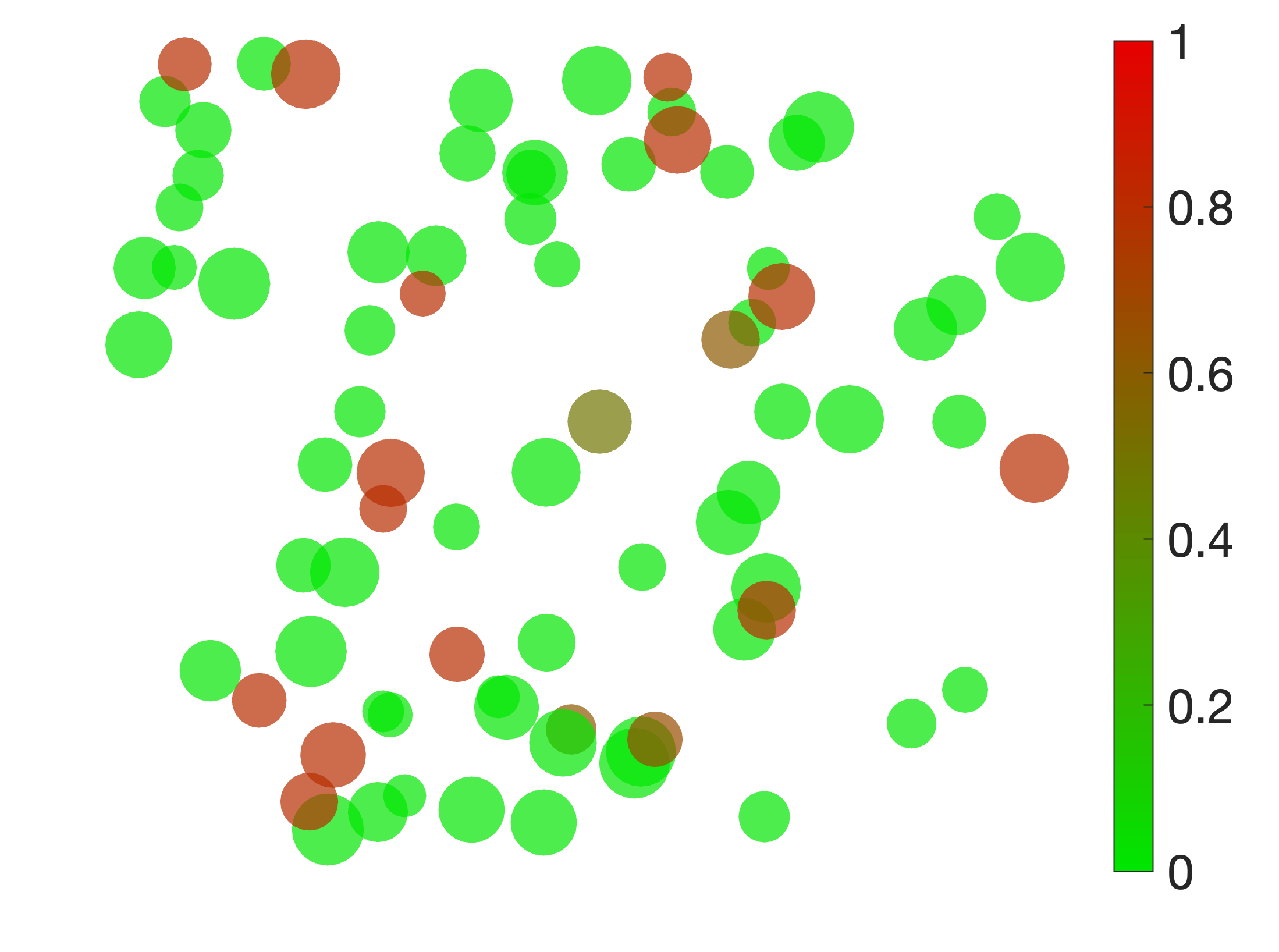}
	     \label{fig:CompareUnmet0.8}
	}  \vspace{-0.2cm}
	\caption{The proportions of unmet demands in different areas}
	\label{fig:CompareUnmet}
	\vspace{-0.2cm}
\end{figure}
Fig.~\ref{fig:CompareUnmet} shows the impact of fairness consideration. 
A lower value of $\beta$ implies that the platform cares more about fairness 
during disruptions. When the fairness condition is relaxed, i.e, larger values of $\beta$ (e.g., $\beta = 0.8$ in Fig.~\ref{fig:CompareUnmet0.8}), there are some areas with particularly high rates of dropped requests (i.e., higher unmet demand). 
On the other hand, when the fairness condition is strictly imposed  (e.g., $\beta=0.2$ in Fig.~\ref{fig:CompareUnmet0.2}), all the areas have comparable proportions of unmet demand.

\begin{table}[ht!]
\centering
\begin{tabular}{|l | c | c|}
\hline
      & $\beta=0.2$            & $\beta=0.8$           \\ \hline
$K=2$ & $\{24,26\}$            & $\{24,26\}$           \\ \hline
$K=4$ & $\{2,11,15,21\}$       & $\{2,10,24,26\}$      \\ \hline
$K=6$ & $\{2,11,15,21,23,26\}$ & $\{2,5,10,24,25,26\}$ \\ \hline
\end{tabular}
\caption{The set of critical ENs}
\vspace{-0.2cm}
    \label{tab:CompareENs}
\end{table}

\begin{figure}[ht!]
    \centering
    \includegraphics[width=0.25\textwidth,height=0.08\textheight]{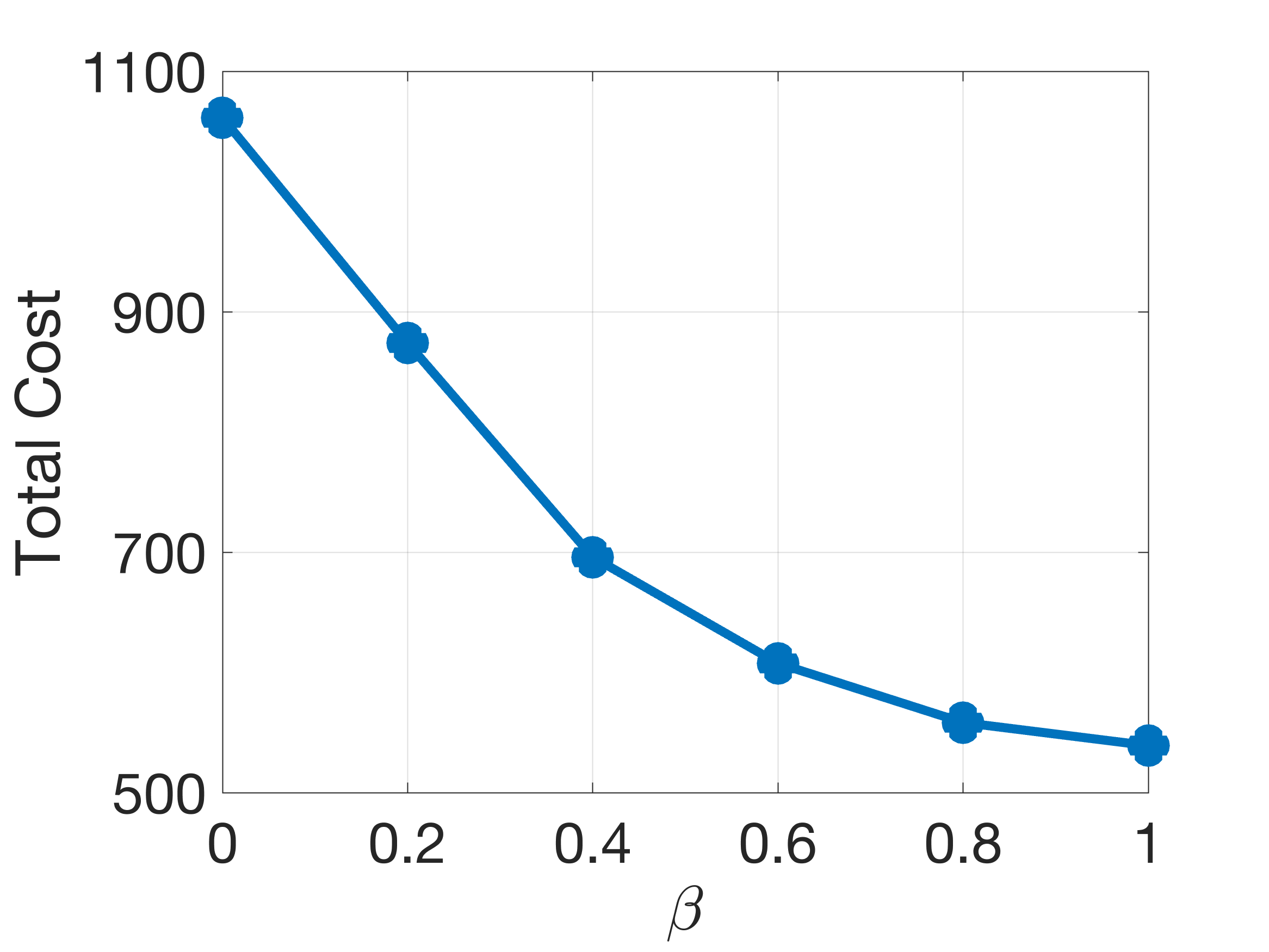}
    \caption{Total penalty for different values of $\beta$}
    \label{fig:CompareCostBeta}
    \vspace{-0.2cm}
\end{figure}

Table \ref{tab:CompareENs} presents the critical EN sets for $\beta$ values of $0.2$ and $0.8$. 
The sets are found to be significantly different, with a varying degree of fairness consideration.
When the fairness constraint is less stringent (i.e., $\beta = 0.8$), the set tends to consist of ENs with high capacity. 
However, when  the fairness constraint is more stringent (i.e., $\beta = 0.2$), 
the set tends to comprise ENs with smaller capacities that have more potential to serve requests from more areas (i.e., close to more areas).
The resource allocation decision during disruption can also differ, even if the sets of critical ENs are the same. 
Indeed, there is an inherent trade-off between fairness and cost, as depicted in Fig.~\ref{fig:CompareCostBeta}. 
Finally, Table~\ref{table:TimeComp} compares the running time of the KKT-based and duality-based reformulations. \textbf{It should be noted that all the experiments were conducted on an M1 MacBook Air with 3.2 GHz processor and 8 GB RAM. On high-performance computing servers, the actual running time can be significantly faster. }

\begin{table}[h]
\centering
\begin{tabular}{|c|c|c|}
\hline
             & KKT (seconds) & Duality (seconds) \\ \hline
$N=10, M=30$ & $15.5208$           & $2.9579$          \\ \hline
$N=30, M=80$ & $235.0888$           & $7.8576$          \\ \hline
$N=80, M=100$ & $4359.9689$           & $27.1822$          \\ \hline
$N=200,~ M=1000$ & N/A & $18599.8247$        \\ \hline
\end{tabular}
\caption{Computational time comparison}\label{table:TimeComp}
\vspace{-0.2cm}
\end{table}

\section{Conclusion}
\label{conc}
This letter introduced a novel fairness-aware robust model to help a budget-constrained EC platform determine the  critical ENs to be safeguarded against possible disruptions. 
The model incorporates the coordination of both hardening and resource allocation decisions in order to enhance system resilience and service quality while striking a balance between fairness and system efficiency.
%
For future research, it is intended to integrate multiple uncertainties, as well as multiple services and resource types into the proposed framework. Exploring the impact of proactive defense through a multi-faceted defender-attacker-defender model is another interesting direction. 

\appendix
\subsection{The Bilevel Attacker-Defender Optimization Problem}\label{app:bilevel}
\vspace{-0.5cm}
\begin{subequations}
\label{eq-Obj:ADfull}
\begin{align}
\textbf{(AD)}~~~~~&  \max_{\bz} ~ (1-\gamma)\sum_{i}\phi_i q_{i} +\gamma\sum_{i,j} d_{i,j} x_{i,j} \nonumber\\
\text{subject to}~~~ & \sum_j z_j\leq K; ~~ z_j \in \{0,~1\}, \quad \forall j \nonumber \\
& \bx, \bq \in \argmin_{\bx, \bq} \Bigg\{ (1-\gamma)\sum_{i}\phi_i q_{i} +\gamma\sum_{i,j} d_{i,j} x_{i,j} \nonumber \\
& \text{subject to} ~ ~~  \sum_{i}x_{i,j}\leq C_j(1-z_j),~\forall j \nonumber \\
&\qquad\qquad\quad  \sum_{j}x_{i,j}+q_{i}= \lambda_i,~ \forall i \nonumber \\
&\qquad\qquad\quad   0 \leq x_{i,j} \leq C_ja_{i,j}   ,~ \forall i,j \nonumber \\
&\qquad\qquad\quad   0 \leq \frac{q_{i}}{\lambda_i} \leq \theta,~ \forall i \nonumber \\
&\qquad\qquad\quad     \Bigg| \frac{q_{i}}{\lambda_i} - \frac{q_{i'}}{\lambda_{i'}} \Bigg| \leq \beta,~ \forall i,~i' \nonumber  \Bigg\}.
\end{align}
\end{subequations}

\subsection{KKT-based Solution}
\label{KKTapproach}

The Lagrangian function of the inner minimization problem 
(\ref{eq:InnerMinimization})-(\ref{eq-const:SP6})
in the subproblem \textbf{(SP)} is:
\begin{subequations} 
\allowdisplaybreaks
\begin{align*} 
\allowdisplaybreaks[2]
\mathcal{L}(\bx,\bq,\pi_j,\mu_i,\sigma_{i,j},\eta_{i,i'},\tau_{i,i'},\nu_{i},\xi_{i,j},\rho_{i}) = -\sum_{i,j}\!\xi_{i,j} x_{i,j} \nonumber\\- \sum_i\! \rho_{i} q_{i} \nonumber
+(1-\gamma)\sum_{i}\phi_i q_{i}+\gamma\sum_{i,j}d_{i,j}x_{i,j}\\ \nonumber
 +\! \sum_j \pi_j \Big(\!\sum_{i}x_{i,j}-\! C_j(1-z_j)\!\Big) \!+ \!\sum_i \mu_i \Big(\sum_{j}x_{i,j}+q_{i} - \lambda_i\!\Big)\\ \nonumber+ \sum_{i,j} \sigma_{i,j} \Big(x_{i,j}- C_ja_{i,j}\Big) + \sum_{i} \nu_{i} \Big(\frac{q_{i}}{\lambda_i} - \theta\!\Big)\cr
    + \sum_{i,i'} \eta_{i,i'} \Big(\frac{q_{i}}{\lambda_i} - \frac{q_{i'}}{\lambda_{i'}} - \beta\Big) - \sum_{i,i'} \tau_{i,i'} \Big(\frac{q_{i}}{\lambda_i} - \frac{q_{i'}}{\lambda_{i'}} + \beta\Big).
\end{align*} 
\end{subequations} 
The KKT conditions give
\begin{subequations}
\begin{align}
\label{eq:kktxijA}
\frac{\partial \mathcal{L}}{\partial x_{i,j}} = \gamma d_{i,j} +\pi_j + \mu_i + \sigma_{i,j} - \xi_{i,j} = 0&,~\forall i,j \\
\label{eq:kktuiA}
\frac{\partial \mathcal{L}}{\partial q_{i}} = (1-\gamma)\phi_i + \mu_{i}  + \frac{1}{\lambda_i}\!\sum_{l=i+1}^M\!\left(\eta_{i,l}-\tau_{i,l}\right)~~&\cr
-\frac{1}{\lambda_i}\sum_{l=1}^{i-1}\!\left(\eta_{l,i}-\tau_{l,i}\right) + \frac{\nu_{i}}{\lambda_i} - \rho_{i} = 0&,~\forall i \\
 \label{eq:kktcs1}
0 \leq \pi_j ~\bot~  C_j(1-z_j) -\sum_{i}x_{i,j} \geq 0&,~\forall j  \\
  \label{eq:kktcs2}
  \sum_{j}x_{i,j}+q_{i} = \lambda_i &, ~\forall j  \\
    \label{eq:kktcs3}
  0 \leq\sigma_{i,j} ~\bot~ C_ja_{i,j} - x_{i,j} \geq 0&, ~\forall j  \\
  \label{eq:kktcs4}
0 \leq \nu_{i} ~\bot~ \theta-\frac{q_{i}}{\lambda_i} \geq 0&,~\forall i,j\\
\label{eq:kktce1}
0 \leq \xi_{i,j} ~\bot~ x_{i,j} \geq 0&,~\forall i,j\\
\label{eq:kktce2}
  0 \leq \rho_{i} ~\bot~ q_{i} \geq 0&,~\forall i
\end{align}
\end{subequations}
where (\ref{eq:kktxijA})-(\ref{eq:kktuiA}) are the stationary conditions, (\ref{eq:kktcs1})-(\ref{eq:kktce2}) are
the primal feasibility, dual feasibility and complementary slackness conditions.
From (\ref{eq:kktxijA})-(\ref{eq:kktuiA}), we have:
\begin{align}
\gamma d_{i,j} +\pi_j + \mu_i + \sigma_{i,j} &= \xi_{i,j}&&,~\forall i,j \label{eq:kktex1}\\
(1-\gamma)\phi_i + \mu_{i}  + \frac{1}{\lambda_i}\sum_{l=i+1}^M\left(\eta_{i,l}-\tau_{i,l}\right)&\cr
-\frac{1}{\lambda_i}\sum_{l=1}^{i-1}\left(\eta_{l,i}-\tau_{l,i}\right) + \frac{\nu_{i}}{\lambda_i} &= \rho_{i} &&,~\forall i.\label{eq:kktex2}
\end{align}
Moreover, from \eqref{eq:kktce1}, if $x_{i,j} > 0$, then $\xi_{i,j}=0$ and \eqref{eq:kktex1} implies that
$\gamma d_{i,j} +\pi_j + \sigma_{i,j} + \mu_i =0,~\forall i,j.$
Thus, 
\beqn
\Big(\gamma d_{i,j} +\pi_j + \sigma_{i,j} + \mu_i\Big)x_{i,j} =0,~\forall i,j.
\eeqn
Similarly, we have 
\beqn
\Bigg[(1-\gamma)\phi_i + \mu_{i}+ \frac{1}{\lambda_i}\sum_{l=i+1}^M\left(\eta_{i,l}-\tau_{i,l}\right)&\cr
-\frac{1}{\lambda_i}\sum_{l=1}^{i-1}\left(\eta_{l,i}-\tau_{l,i}\right) + \frac{\nu_{i}}{\lambda_i}\Bigg]q_{i} &=0,~\forall i.
\eeqn
Therefore, the conditions \eqref{eq:kktce1}-\eqref{eq:kktex2} are equivalent to the following set of constraints:
\begin{subequations}
\begin{align}
0 \leq \gamma d_{i,j} + \pi_j + \mu_i + \sigma_{i,j}  ~\bot~ x_{i,j} \geq 0&,~\forall i,j \label{eq:kktxijex} \\
\label{eq:kktuiex}
0 \leq (1-\gamma)\phi_i + \mu_{i} + \frac{1}{\lambda_i}\sum_{l=i+1}^M\left(\eta_{i,l}-\tau_{i,l}\right)~~~~&\cr
-\frac{1}{\lambda_i}\sum_{l=1}^{i-1}\left(\eta_{l,i}-\tau_{l,i}\right) + \frac{\nu_{i}}{\lambda_i} ~\bot~ q_{i} \geq 0 &,~\forall i.
\end{align}
\end{subequations}
It is worth noting that we can directly use the  set of KKT conditions  (\ref{eq:kktex1})-(\ref{eq:kktex2})
 to solve the subproblem, but it will involve more variables (i.e., $\xi$ and $\rho$) compared to solving the subproblem with (\ref{eq:kktxijex})-(\ref{eq:kktuiex}). In brief, based on  the Karush–Kuhn–Tucker (KKT) conditions above, the (\textbf{AD}) problem (\ref{ADProb}) is equivalent to the following problem with complementary constraints:
\begin{subequations}
\begin{align}
\label{eq:SPobjf}
\max_{\bx,\bq,\bz,\pi,\sigma,\mu,\eta,\tau,\nu}  &(1-\gamma)\sum_{i}\phi_i q_{i} + \gamma\sum_{i,j}d_{i,j}x_{i,j}\\
\text{subject to}~~~
    &\eqref{eq:kktcs1}, \eqref{eq:kktcs3}, \eqref{eq:kktcs2},\eqref{eq:kktcs4}, \eqref{eq:kktxijex},\eqref{eq:kktuiex} \label{eq:AllCSConst} \\ 
    & z_j\in\{0,1\},~\forall j,~~\sum_jz_j\leq K, \label{eq:AttackConst}
\end{align}
\end{subequations}
where the last constraints in \eqref{eq:AttackConst} represents the set of feasible attack plans $\mathcal{Z}$. 

\textbf{The Complete MILP Formulation}
Note that a complimentary constraint $0 \leq x \bot \pi \geq 0$ means $x \geq 0, \pi \geq$ and $x . \pi = 0$. Thus, it is a nonlinear constraint. However, this nonlinear complimentary constraint can be transformed into equivalent exact linear constraints by using the Fortuny-Amat transformation. 
Specifically, the complementarity  condition $0 \leq x \bot \pi \geq 0$ is  equivalent  to the following set of mixed-integer linear constraints:
\begin{subequations}
\begin{align}
x \geq 0;~~ x \leq (1-u)M \\
\pi \geq 0;~~ \pi \leq  uM,~~ u\in\{0;1\},
\end{align}
\end{subequations}
where M is a sufficiently large constant.
By applying this transformation to all the complementary constraints listed in (\ref{eq:AllCSConst}), we obtain an MILP that is equivalent to the problem (\ref{ADProb}). The explicit form of this  MILP is as follows:
\begin{subequations} 
\allowdisplaybreaks
\begin{align}
\label{eq-Obj:SP-MILP}
&\max_{\bx,\bq,\bz,\pi,\sigma,\mu,\eta,\tau,\nu,\bu}  (1-\gamma)\sum_{i}\phi_i q_{i} ~+&& \!\!\!\!\!\!\!\!\!\!\!\!\!\!\gamma\sum_{i,j}d_{i,j}x_{i,j}   \\
&\text{subject to}&&~\nonumber\\
&0 \leq  \gamma d_{i,j} + \pi_j + \mu_i + \sigma_{i,j} \leq u_{i,j}^0 M_{i,j}^0 &&,\forall i,j\\
&0 \leq x_{i,j} \leq (1-u_{i,j}^0) M_{i,j}^0 &&,\forall i,j\\
&0 \leq(1-\gamma)\phi_i + \mu_{i} + \frac{1}{\lambda_i}\sum_{l=i+1}^M\left(\eta_{i,l}-\tau_{i,l}\right)&&\cr
&~~~~-\frac{1}{\lambda_i}\sum_{l=1}^{i-1}\left(\eta_{l,i}-\tau_{l,i}\right) + \frac{\nu_{i}}{\lambda_i} \leq u_{i}^1 M_i^1 &&,\forall i\\
&0 \leq q_{i} \leq (1-u_{i}^1) M_i^1 &&,\forall i\\
&0 \leq C_j(1-z_j) -\sum_{i}x_{i,j} \leq u_j^2 M_j^2 &&,\forall j\\
&0 \leq  \pi_j \leq (1-u_j^2) M_j^2 &&,\forall j\\ 
&0 \leq C_ja_{i,j} - x_{i,j} \leq u_{i,j}^3 M_{i,j}^3 &&,\forall i,j  \\ 
&0 \leq \sigma_{i,j} \leq (1-u_{i,j}^3) M_{i,j}^3 &&,\forall i,j  \\
&0 \leq \beta-\frac{q_{i}}{\lambda_i} + \frac{q_{i'}}{\lambda_{i'}} \leq u_{i,i'}^4 M_{i,i'}^4 &&,\forall i,i'  \\ 
&0 \leq \eta_{i,i'} \leq (1-u_{i,i'}^4) M_{i,i'}^4 &&,\forall i,i'  \\
&0 \leq \frac{q_{i}}{\lambda_i} - \frac{q_{i'}}{\lambda_{i'}}+\beta \leq u_{i,i'}^5 M_{i,i'}^5 &&,\forall i,i'  \\ 
&0 \leq \tau_{i,i'} \leq (1-u_{i,i'}^5) M_{i,i'}^5 &&,\forall i,i'  \\
&0 \leq \theta-\frac{q_{i}}{\lambda_i} \leq u_{i}^6 M_{i}^6 &&,\forall i  \\ 
&0 \leq \nu_{i} \leq (1-u_{i}^6) M_{i}^6 &&,\forall i  \\
&    \sum_{j}x_{i,j}+q_{i} = \lambda_i  &&,\forall i\\
&u_{i,j}^0, u_{i}^1, u_{j}^2, u_{i,j}^3,u_{i,i'}^4,u_{i,i'}^5,u_{i}^6 \in \{0,1\} &&,\forall i,i',j\\
&\sum_{j}z_j\leq K, ~~z_j\in\{0,1\}&&,\forall j
\label{MILPSPeq3}
\end{align}
\end{subequations} 
where $\bu$ represents the set of binary variables $u^0, u^1, u^2$, $u^3, u^4, u^5, u^6$. Also,  $M_{i,j}^0, M_{i}^1, M_j^2, M_{i,j}^3, M_{i,i'}^4,M_{i,i'}^5 , M_{i}^6$ are sufficiently large numbers. The value of each $M$ should be large enough to ensure feasibility of the associated constraint. On the other hand, the value of each $M$ should not be too large to enhance the computational speed of the solver. Indeed, the value of each $M$ should be tighten to the limits of parameters and variables in the corresponding constraint.

\subsection{Flowchart of the Proposed Solution Approach} 
\label{app:Flowchart}
\begin{figure}[ht!]
    \centering
    \includegraphics[width=0.5\textwidth]{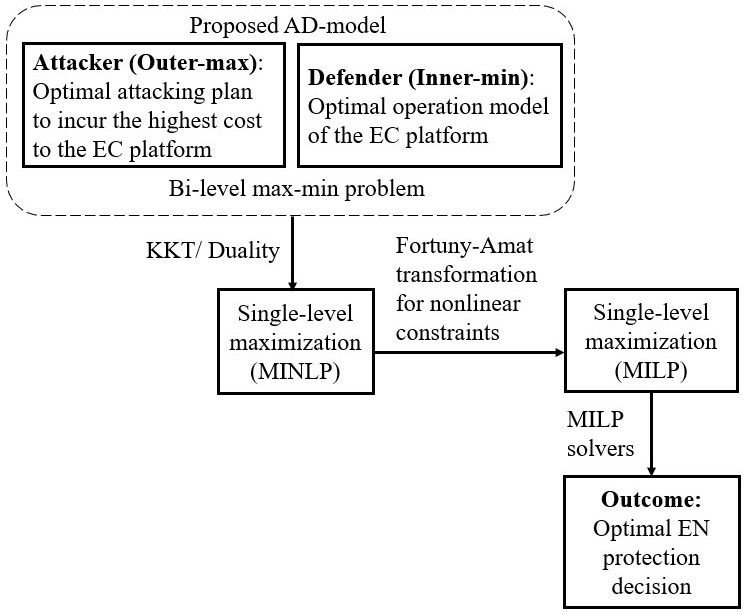}
    \caption{Flowchart of the solution approach}
 	\label{fig:Flowchart}
\end{figure}

The flow chart shown in Figure~\ref{fig:Flowchart} summarizes the proposed system model and solution approach. Specifically, the inner minimization problem in (3) is first converted as a linear minimization program shown in (6). Then, by using\textit{ strong duality in linear programming}, (6) is equivalent to the non-linear maximization problem (7). At this step, the bilevel maxmin problem becomes a single-level non-linear maximization problem. The Fortuny-Amat transformation (big-M method) is then employed to convert bilinear  constraints into  equivalent sets of linear equations. Consequently, the initial bilevel problem (3) is transformed into a mixed-integer linear program (MILP) as shown in (9), which can be solved efficiently by off-the-shelf MILP solvers such as Gurobi, Mosek, Cplex. Note that instead of using LP strong duality theorem, we can also use KKT conditions to convert (6) into an equivalent set of linear equations. The KKT-based solution is presented in \textit{Appendix B}. However, as explained, the KKT-based approach results in more integer variables and constraints in the MILP reformulation compared to that of the duality-based approach. The solution to the reformulated MILP problem is an optimal set of critical ENs that should be protected.

\end{document}